\documentclass[11pt]{article}
\pdfoutput=1 
\usepackage{amsfonts, amsmath, amssymb, amsthm, fullpage, enumitem, mathtools, microtype, CJK, hyperref}
\usepackage[noabbrev, capitalize]{cleveref}

\title{On symmetric hollow integer matrices\\with eigenvalues bounded from below}
\author{
    Zilin Jiang \begin{CJK}{UTF8}{gkai}姜子麟\end{CJK}\thanks{School of Mathematical and Statistical Sciences, and School of Computing and Augmented Intelligence, Arizona State University, Tempe, AZ 85281, USA. Email: {\tt zilinj@asu.edu}.}
}
\date{}

\newtheorem{theorem}{Theorem}
\newtheorem{lemma}[theorem]{Lemma}
\newtheorem{proposition}[theorem]{Proposition}
\newtheorem{problem}[theorem]{Problem}

\theoremstyle{definition}
\newtheorem{definition}[theorem]{Definition}

\theoremstyle{remark}
\newtheorem*{remark}{Remark}

\DeclarePairedDelimiter\abs{\lvert}{\rvert}%

\newcommand{\dset}[2]{\left\{{#1}\colon{#2}\right\}}
\newcommand{\sset}[1]{\left\{{#1}\right\}}

\newcommand{\al}{\alpha}
\newcommand{\be}{\beta}
\newcommand{\la}{\lambda}
\newcommand{\las}{\lambda^*}

\newcommand{\Sla}{\mathcal{S}(\lambda)}
\newcommand{\Stwo}{\mathcal{S}(2)}
\newcommand{\Mla}{\mathcal{M}(\lambda)}
\newcommand{\M}{\mathcal{M}}
\newcommand{\N}{\mathbb{N}}

\begin{document}

\maketitle

\begin{abstract}
    A hollow matrix is a square matrix whose diagonal entries are all equal to zero. Define $\lambda^* = \rho^{1/2} + \rho^{-1/2} \approx 2.01980$, where $\rho$ is the unique real root of $x^3 = x + 1$. We show that for every $\lambda < \lambda^*$, there exists $n \in \mathbb{N}$ such that if a symmetric hollow integer matrix has an eigenvalue less than $-\lambda$, then one of its principal submatrices of order at most $n$ does as well. However, the same conclusion does not hold for any $\lambda \ge \lambda^*$.
\end{abstract}

\noindent\textbf{Keywords:} Symmetric hollow integer matrix; Minimal forbidden submatrix

\noindent\textbf{Mathematics Subject Classification:} 05C50; 15A18

\section{Introduction} \label{sec:intro}

The Cauchy interlacing theorem implies that if the eigenvalues of a Hermitian matrix $A$ are at least $-2$, then the same is true for every principal submatrix of $A$, particularly for those of order up to, say, $10$.

In \cite{V}, Vijayakumar proved the following surprising converse for certain matrices: If a symmetric hollow integer matrix has an eigenvalue less than $-2$, then one of its principal submatrices of order at most $10$ does as well. To state this result in a different way, we introduce the following concepts.

\begin{definition}[Minimal forbidden submatrix]
    Given a family $\M$ of symmetric hollow integer matrices, a symmetric hollow integer matrix $A$ is a \emph{minimal forbidden submatrix} for $\M$ if $A$ itself is not in $\M$ but every proper principal submatrix of $A$ is in $\M$. Let $\Mla$ be the family of symmetric hollow integer matrix with smallest eigenvalue at least $-\la$.
\end{definition}

The aforementioned result amounts to saying that the order of every minimal forbidden submatrix for $\M(2)$ is at most $10$.

The key observation in \cite{V} is that a minimal forbidden submatrix $A$ for $\M(2)$ of order more than $3$ cannot have any entry bigger than $1$ in absolute value. Thus $A$ can be viewed as the signed adjacency matrix of a signed graph, say $F$, and moreover $F$ is a minimal forbidden subgraph for $\Stwo$. Here, minimal forbidden subgraph and the family $\Sla$ of signed graphs are defined analogously.

\begin{definition}[Minimal forbidden subgraph]
    Given a family $\mathcal{S}$ of signed graphs, a signed graph $F$ is a \emph{minimal forbidden subgraph} for $\mathcal{S}$ if $F$ itself is not in $\mathcal{S}$ but every proper induced subgraph of $F$ is in $\mathcal{S}$. Let $\Sla$ be the family of signed graphs whose smallest eigenvalue of the signed adjacency matrix is at least $-\la$
\end{definition}

Indeed, it is proved in \cite{V} that every minimal forbidden subgraph for the family $\Stwo$ has at most $10$ vertices.

Inspired by the connection between symmetric hollow integer matrices and signed graphs, we capitalize on the recent developments on signed graphs in \cite{JP} to qualitatively generalize Vijayakumar's matrix result as follows.

\begin{theorem} \label{thm:main}
    For every $\la < \las$, there exists $n \in \N$ such that if a symmetric hollow integer matrix has an eigenvalue less than $-\la$, then one of its principal submatrices of order at most $n$ does as well. However, the same conclusion does not hold for any $\la \ge \las$.
\end{theorem}

The main difficulty occurs when $\la \in (2, \las)$ --- a minimal forbidden submatrix for $\Mla$ could have entries bigger than $1$ in absolute value. We overcome this difficulty by characterizing these minimal forbidden submatrices, and bounding their order using Dickson's lemma.

\begin{lemma}[Lemma~A of Dickson~\cite{D}] \label{lem:dickson}
    The partially ordered set $(\N^n, \le)$, in which $(a_1, \dots, a_n) \le (b_1, \dots, b_n)$ if and only if $a_i \le b_i$ for every $i$, does not contain infinite antichains. \qed
\end{lemma}

All of the proofs will be given in \cref{sec:proof}. We conclude in \cref{sec:rem} with some brief remarks.

\section{Proofs} \label{sec:proof}

\subsection{Proof for $\la \le 2$}

A useful tool in spectral graph theory for signed graphs is switching --- two signed graphs are \emph{switching equivalent} if one graph can be obtained from the other by reversing all the edges in a cut-set. We extend this notion to matrices, and we make the simple observation on minimal forbidden submatrices that are switching equivalent.

\begin{definition}[Switching equivalence]
    Two matrices $A$ and $B$ are \emph{switching equivalent} if there exists a diagonal matrix $D$ with $\pm 1$ diagonal entries such that $A = D^\intercal BD$.
\end{definition}

\begin{proposition} \label{prop:switching}
    Suppose that two matrices $A$ and $B$ are switching equivalent. For every $\la \in \mathbb{R}$, the matrix $A$ is a minimal forbidden submatrix of $\M(\la)$ if and only if the matrix $B$ is as well. \qed
\end{proposition}

Recently, motivated by a discrete-geometric question on spherical two-distance sets in \cite{JJ}, Polyanskii and the author determined the set of $\la$ for which there are finitely many minimal forbidden subgraphs for $\Sla$.

\begin{theorem}[Theorem~1.5 of Jiang and Polyanskii~\cite{JP}] \label{thm:forb}
    The family $\Sla$ of signed graphs with smallest eigenvalue at least $-\la$ has a finite set of minimal forbidden subgraphs if and only if $\la < \las$. \qed
\end{theorem}

We need one last ingredient for the case where $\la \le 2$. One can verify by brute force the following fact, which already appeared in \cite[p.~111]{V}.

\begin{proposition} \label{thm:root5}
    For all $a, b \in \sset{0,\pm1, \pm2}$, if $a \neq b$, then the smallest eigenvalue of
    \[
        \begin{pmatrix}
            0 & 2 & a \\
            2 & 0 & b \\
            a & b & 0
        \end{pmatrix}
    \]
    is at most $-\sqrt{5}$. \qed
\end{proposition}

\begin{proof}[Proof of \cref{thm:main} for $\la \le 2$]
    Let $\mathcal{F}$ be the set of minimal forbidden subgraphs for $\Sla$. With hindsight, set $n = \max \left(\sset{3} \cup \dset{\abs{F}}{F \in \mathcal{F}}\right)$. Note that $n$ is finite because $\mathcal{F}$ is finite according to \cref{thm:forb}. Suppose that $A = (a_{ij})$ is a minimal forbidden submatrix for $\M(\la)$. We shall prove that the order of $A$ is at most $n$.

    Suppose for a moment that $\abs{a_{ij}} > \la$ for some $i < j$. Then the eigenvalues of the principal submatrix
    \[
        \begin{pmatrix}
            a_{ii} & a_{ij} \\
            a_{ji} & a_{jj}
        \end{pmatrix}
    \]
    are $\pm a_{ij}$, one of which is less than $-\la$. Hereafter, we may assume that every entry of $A$ does not exceed $\la$ in absolute value. We break the rest of the proof into two cases.

    \paragraph{Case 1:} $\abs{a_{ij}} = 2$ for some $i < j$. In this case, $\la = 2$. Since $A$ is switching equivalent to a matrix where $a_{ij} = 2$ through the diagonal matrix where all diagonal entries are $+1$ except the $i$-th is $-1$, in view of \cref{prop:switching}, we may assume that $a_{ij} = 2$.
    
    Without loss of generality, we may further assume that $(i,j) = (1,2)$ or equivalently $a_{12} = 2$, and the order of $A$ is more than $3$. Under these assumptions, according to \cref{thm:root5}, we have that $a_{1k} = a_{2k}$ for every $k \ge 3$ because otherwise we would have found a principal submatrix of order $3$ that has an eigenvalue less than $-2$.
    
    Let $A'$ be the principal submatrix obtained by removing the first row and the first column of $A$. Since $A$ is a minimal forbidden submatrix for $\M(2)$, the matrix $A'$ is in $\M(2)$, and so $A' + 2I$ is positive semidefinite. Since the first two rows and the first two columns of $A + 2I$ are the same, one can check that $A + 2I = B^\intercal (A' + 2I) B$, where 
    \[
        B = \begin{pmatrix}
            1 & 1 & & & & \\
            & & 1 & & & \\
            & & & \ddots & \\
            & & & & 1
        \end{pmatrix}
    \]
    and so $A + 2I$ is positive semidefinite, which contradicts with the assumption that $A \not\in \M(2)$.

    \paragraph{Case 2:} $\abs{a_{ij}} \le 1$ for every $i < j$. Then $A$ is the signed adjacency matrix of a signed graph, say $F$, and moreover $F \in \mathcal{F}$. Clearly, the order of $A$, that is $\abs{F}$, is at most $n$.
\end{proof}

\begin{remark}
    For a long period of time, the author of the current paper had trouble understanding the second paragraph in the proof of Theorem~4.2 in \cite{V}. Our treatment of the case where $\abs{a_{ij}} = 2$ and $\la = 2$ in the proof of \cref{thm:main} clarifies that argument in \cite{V}.
\end{remark}

\subsection{Proof for $\la \in (2, \las)$}

For the case where $\la \in (2, \las)$, we characterize the minimal forbidden submatrices for $\Mla$. To that end, we introduce the following symmetric hollow integer matrix.

\begin{definition}
    Given a signed graph $F$ and a mapping $a \colon V(F) \to \N^+$, the matrix $A_{F,a}$ is indexed by the set $\dset{(u,i)}{u \in V(F) \text{ and } i \in \sset{1, \dots, a(u)}}$, and its entries are defined by
    \[
        A_{F, a}((u,i), (v,j)) = \begin{cases}
            0 & \text{if }u = v \text{ and } i = j; \\
            +2 & \text{if }u = v \text{ but } i \neq j; \\
            +1 & \text{if }u \neq v \text{ and } uv \text{ is a positive edge of }F; \\
            -1 & \text{if }u \neq v \text{ and } uv \text{ is a negative edge of }F. \\
        \end{cases}
    \]
\end{definition}

\begin{lemma} \label{lem:char}
    For every symmetric hollow integer matrix $A$ of order at least $3$, if each principal submatrix of $A$ of order $3$ has its smallest eigenvalue more than $-\sqrt5$, then there exist a signed graph $F$ and $a \colon V(F) \to \N^+$ such that $A$ is switching equivalent to $A_{F,a}$.
\end{lemma}

\begin{proof}
    Suppose for a moment that $\abs{a_{ij}} \ge 3$ for some $i < j$. Then the eigenvalues of the principal submatrix
    \[
        \begin{pmatrix}
            a_{ii} & a_{ij} \\
            a_{ji} & a_{jj}
        \end{pmatrix}
    \]
    are $\pm a_{ij}$, one of which is less than $-\sqrt{5}$. Thus any principal submatrix of order $3$ containing the above principal submatrix of order $2$ has an eigenvalue less than $-\sqrt{5}$, which is a contradiction. Hereafter, we may assume that every entry of $A$ does not exceed $2$ in absolute value.

    Suppose that the matrix $A$ is indexed by $V$. Let $G$ be the graph on $V$ where $ij$ is an edge if and only if $\abs{a_{ij}} = 2$. Let the connected components of $G$ partition $V$ into disjoint sets $V_1, \dots, V_m$.

    Consider each vertex subset $V_k$. Since the induced subgraph $G[V_k]$ is connected, let $T_k$ be a spanning tree of $G[V_k]$. Since $T_k$ is a tree, there exists a partition of $V_k$ into $X_k$ and $Y_k$ such that for every edge $ij \in T_k$, $a_{ij} = -2$ if and only if $i \in X_k$ and $j \in Y_k$.

    Now, define the diagonal matrix $D$ by
    \[
        d_{ii} = \begin{cases}
            +1 \text{ if } i \in X_1 \cup \dots \cup X_m; \\
            -1 \text{ if } i \in Y_1 \cup \dots \cup Y_m, \\
        \end{cases}
    \]
    and define the matrix $B = (b_{ij}) = D^\intercal A D$, which is switching equivalent to $A$. Due to the switching, $b_{ij} = 2$ for every edge $ij \in T_k$. Through repeated application of \cref{thm:root5}, one can show that
    \begin{enumerate}[label=(\alph*)]
        \item for every $i, j \in V_k$ with $i \neq j$, the entry $b_{ij} = 2$, and
        \item for every $i \in V_k$ and $j \in V_\ell$ with $k \neq \ell$, the entry $b_{ij}$ is a constant, say $c_{k\ell} \in \sset{-1,0,+1}$, that depends only on $k$ and $\ell$.
    \end{enumerate}

    Finally, one can check that $B = A_{F,a}$, where $F$ is the signed graph on $\sset{1,\dots, m}$ whose signed adjacency matrix is $(c_{k\ell})$, and $a\colon \sset{1,\dots, m} \to \N^+$ is defined by $a(k) = \abs{V_k}$.
\end{proof}

In addition to \cref{thm:forb}, we need the following finiteness result.

\begin{theorem}[Theorem 3.9 of Jiang and Polyanskii~\cite{JP}] \label{thm:finite}
    For every $\la \in [2, \las)$, the number of connected signed graphs in $\Sla \setminus \Stwo$ is finite. \qed
\end{theorem}

We are now in the position to prove for $\la \in (2,\las)$.

\begin{proof}[Proof of \cref{thm:main} for $\la \in (2, \las)$]
    Let $\mathcal{F}$ be the set of minimal forbidden subgraphs for $\Sla$. With hindsight, set
    \[
        n = \max\left(\sset{3} \cup \dset{\abs{F}}{F \in \mathcal{F}} \cup \dset{n_F}{F \text{ is a connected signed graph in } \Sla \setminus \Stwo}\right),
    \]
    where $n_F$ is the maximum order of a minimal forbidden submatrix for $\Mla$ of the form $A_{F, a}$.
    
    To justify that $n$ is finite, because both $\mathcal{F}$ and the set of connected signed graphs in $\Sla \setminus \Stwo$ are finite according to \cref{thm:forb,thm:finite}, it is enough to show that $n_F$ is finite for every signed graph $F$. Observe that if $a \le b$, then $A_{F, a}$ is a principal submatrix of $A_{F, b}$. Since no minimal forbidden submatrix for $\Mla$ is a principal submatrix of another, for a fixed signed graph $F$, the set
    \[
        \dset{a \colon V(F) \to \N^+}{A_{F, a} \text{ is a minimal forbidden submatrix for }\Mla}
    \]
    forms an antichain. According to \cref{lem:dickson}, this antichain is finite, and so is $n_F$.

    Suppose that $A = (a_{ij})$ is a minimal forbidden submatrix for $\M(\la)$. We shall prove that the order of $A$ is at most $n$. Without loss of generality, we may assume that the order of $A$ is more than $3$. Since $\la < \las < \sqrt{5}$, each principal submatrix of $A$ of order $3$ has its smallest eigenvalue more than $-\sqrt{5}$. According to \cref{lem:char}, there exist a signed graph $F$ and $a \colon V(F) \to \N^+$ such that $A$ is switching equivalent to $A_{F, a}$. We break the rest of the proof into two cases.

    \paragraph{Case 1:} $\abs{a_{ij}} = 2$ for some $i < j$. In view of \cref{prop:switching}, we may assume that $A = A_{F, a}$. Since $A$ is a minimal forbidden submatrix for $\Mla$, the matrix $A$ must be irreducible, which implies that the signed graph $F$ is connected. Since $\abs{a_{ij}} = 2$, we know that the signed adjacency matrix $A_F$ of the signed graph $F$ is a proper principal submatrix of $A_{F, a}$, and so $F \in \Sla$.

    We claim that $F \not\in \Stwo$. Assume for the sake of contradiction that $A_F + 2I$ is positive semidefinite.
    Let $B$ be the matrix whose rows are indexed by $V(F)$, whose columns are indexed by $\dset{(v, i)}{v \in V(F) \text{ and }i \in \sset{1, \dots, a(u)}}$, and whose entries are defined by
    \[
        B(u, (v, i)) = \begin{cases}
            1 & \text{if }u = v;\\
            0 & \text{otherwise}.
        \end{cases}
    \]
    One can check that $A_{F, a} + 2I = B^\intercal (A_F + 2I) B$, and so $A_{F,a} + 2I$ is positive semidefinite, which contradicts with the assumptions that $A \not\in \Mla$ and $\la > 2$.

    Therefore $F$ is a connected signed graph in $\Sla \setminus \Stwo$, and so the order of the minimal forbidden matrix $A_{F,a}$ for $\Mla$ is at most $n_F$, which is at most $n$.

    \paragraph{Case 2:} $\abs{a_{ij}} \le 1$ for every $i < j$. Then $A$ is the signed adjacency matrix of the signed graph $F$, and moreover $F \in \mathcal{F}$. Clearly, the order of $A$, that is $\abs{F}$, is at most $n$.
\end{proof}

\subsection{Proof for $\la \ge \las$}

\begin{proof}[Proof of \cref{thm:main} for $\la \ge \las$]
    Because of \cref{thm:forb}, there exists a sequence of minimal forbidden subgraphs $F_1, F_2, \dots$ for $\Sla$ with increasing number of vertices. Let $A_n$ be the signed adjacency matrix of $F_n$. By the choice of $F_n$, we know that $A_n$ is a minimal forbidden submatrix for $\Mla$.
\end{proof}

\section{Concluding remarks} \label{sec:rem}

We have determined the set of $\la$ for which the maximum order of a minimal forbidden submatrix for $\Mla$ is finite. It would be interesting to investigate other natural matrix classes.

We highlight the following matrix class. Consider the family $\M'(\la)$ of symmetric integer matrices with spectral radius at most $\la$.

\begin{problem} \label{prob:matrix}
    Determine the set of $\la$ for which the maximum order of a minimal forbidden submatrix for $\M'(\la)$ is finite.
\end{problem}

In this direction, motivated by the Lehmer's Mahler measure problem, McKee and Smyth enumerated in \cite{MS} the minimal forbidden submatrices for $\M'(2)$, all of which are of order at most $10$. 

A good starting point to attack \cref{prob:matrix} is to understand the minimal forbidden subgraphs for the family $\mathcal{S}'(\la)$ of signed graphs with spectral radius at most $\la$.

\begin{problem} \label{prob:signed-graph}
    Determine the set of $\la$ for which $\mathcal{S}'(\la)$ has a finite set of minimal forbidden subgraphs.
\end{problem}

For unsigned graphs, such a set of $\la$ has been determined in \cite{JP20} --- the family of graphs with spectral radius at most $\la$ has a finite set of minimal forbidden subgraphs if and only if $\la < \la'$ and $\la \not\in \sset{\al_2, \al_3, \dots}$, where $$\la' = \sqrt{2+\sqrt{5}} \approx 2.05817, \quad \al_m = \be_m^{1/2} + \be_m^{-1/2},$$ and $\be_m$ is the largest root of $x^{m+1}(x-1) = x^m-1$.

The above result on unsigned graphs implies that the answer to \cref{prob:signed-graph} must satisfy that $\la < \la'$ and $\la \not\in \sset{\al_2, \al_3, \dots}$. Another potentially useful result is due to Wang, Dong, Hou and Li, who classified in \cite{WDHL} all the connected signed graphs with spectral radius at most $\la'$.

Finally, we reiterate the following problem initially raised by McKee and Smyth \cite[p.~106]{MS}.

\begin{problem}
    Classify all the irreducible symmetric integer matrices with spectral radius in $(2, \la')$.
\end{problem}

\section*{Acknowledgements}

The author would like to thank an anonymous referee who read the first arXiv version of \cite{JP} and made the author notice the gap in the proof of \cite[Theorem~4.2]{V}. The author is also grateful to Theodore Gossett and Shengtong Zhang who read the preliminary version of the paper and pointed out many inaccuracies.

\bibliographystyle{plain}
\bibliography{matrix}

\begin{thebibliography}{1}

\bibitem{D}
Leonard~Eugene Dickson.
\newblock Finiteness of the odd perfect and primitive abundant numbers with {$n$} distinct prime factors.
\newblock {\em Amer. J. Math.}, 35(4):413--422, 1913.

\bibitem{JP20}
Zilin Jiang and Alexandr Polyanskii.
\newblock Forbidden subgraphs for graphs of bounded spectral radius, with applications to equiangular lines.
\newblock {\em Israel J. Math.}, 236(1):393--421, 2020.
\newblock \href{http://arxiv.org/abs/1708.02317}{\tt arXiv:1708.02317 [math.CO]}.

\bibitem{JP}
Zilin Jiang and Alexandr Polyanskii.
\newblock Forbidden induced subgraphs for graphs and signed graphs with eigenvalues bounded from below, 2024.
\newblock \href{http://arxiv.org/abs/2111.10366}{\tt arXiv:2111.10366 [math.CO]}.

\bibitem{JJ}
Zilin Jiang, Jonathan Tidor, Yuan Yao, Shengtong Zhang, and Yufei Zhao.
\newblock Spherical two-distance sets and eigenvalues of signed graphs.
\newblock {\em Combinatorica}, 43(2):203--232, 2023.
\newblock \href{http://arxiv.org/abs/2006.06633}{\tt arXiv:2006.06633 [math.CO]}.

\bibitem{MS}
James McKee and Chris Smyth.
\newblock Integer symmetric matrices of small spectral radius and small {M}ahler measure.
\newblock {\em Int. Math. Res. Not. IMRN}, (1):102--136, 2012.
\newblock \href{http://arxiv.org/abs/0907.0371}{\tt arXiv:0907.0371 [math.NT]}.

\bibitem{V}
G.~R. Vijayakumar.
\newblock Signed graphs represented by {$D_\infty$}.
\newblock {\em European J. Combin.}, 8(1):103--112, 1987.

\bibitem{WDHL}
Dijian Wang, Wenkuan Dong, Yaoping Hou, and Deqiong Li.
\newblock On signed graphs whose spectral radius does not exceed $\sqrt{2+\sqrt 5}$.
\newblock {\em Discret. Math.}, 346:113358, 2022.
\newblock \href{http://arxiv.org/abs/2203.01530}{\tt arXiv:2203.01530 [math.CO]}.

\end{thebibliography}

\end{document}